\documentclass[12pt, a4paper]{amsart}
\usepackage{amssymb, amsthm, graphicx, epsfig, verbatim}
\newtheorem{thm}{Theorem}[section]

\newtheorem{lem}[thm]{Lemma}

\theoremstyle{definition}

\newcommand{\al}{\alpha}

\newcommand{\Om}{\Omega}
\newcommand{\om}{\omega}

\newcommand{\si}{\sigma}
\newcommand{\Si}{\Sigma}

\newcommand{\ze}{\zeta}

\newcommand{\x}{\times}
\newcommand{\s}{\mathbf s}

\renewcommand{\t}{\mathbf t}

\newcommand{\C}{\mathbb C}
\newcommand{\Z}{\mathbb Z}
\newcommand{\N}{\mathbb N}
\newcommand{\Q}{\mathbb Q}

\newcommand{\CP}{{\mathbb C}{\mathbb P}}

\newcommand{\del}{\partial}

\newcommand{\co}{\colon\thinspace}

\begin{document}
\mathsurround=1pt
\title{Tight, not semi--fillable contact circle bundles}

\author{Paolo Lisca}
\address{Dipartimento di Matematica\\
Universit\`a di Pisa \\I-56127 Pisa, Italy}

\author{Andr\'{a}s I. Stipsicz}
\address{R\'enyi Institute of Mathematics\\
Hungarian Academy of Sciences\\
H-1053 Budapest\\
Re\'altanoda utca 13--15, Hungary}
\email{}

\begin{abstract}
Extending our earlier results, we prove that certain tight contact
structures on circle bundles over surfaces are not symplectically
semi--fillable, thus confirming a conjecture of Ko Honda.
\end{abstract}
\thanks{E-mail addresses: lisca@dm.unipi.it (P.~Lisca), 
stipsicz@math-inst.hu (A.I.~Stipsicz)\newline
The first author was partially supported by
  MURST, and he is a member of EDGE, Research Training Network
  HPRN-CT-2000-00101, supported by The European Human Potential
  Programme. 
  The second author was partially supported by Sz\'echenyi
  Professzori \"Oszt\"ond{\'\i}j and OTKA}

\maketitle

\section{Introduction}\label{s:intro}
Let $Y$ be a closed, oriented three--manifold. A \emph{positive,
coorientable contact structure} on $Y$ is the kernel
$\xi=\ker\al\subset TY$ of a one--form $\alpha\in\Om^1(Y)$ such that
$\alpha\wedge d\alpha$ is a positive volume form on $Y$. The pair
$(Y,\xi)$ is a~\emph{contact three--manifold}. In this paper we only
consider positive, coorientable contact structures, so we call them
simply `contact structures'. For an introduction to contact structures
the reader is referred to~\cite{Ae}, Chapter~8 and~\cite{E}.

There are two kinds of contact structures $\xi$ on $Y$. If there
exists an embedded disc $D\subset Y$ tangent to $\xi$ along its
boundary, $\xi$ is called~\emph{overtwisted}, otherwise it is said to
be \emph{tight}. The isotopy classification of overtwisted contact
structures coincides with their homotopy classification as tangent
two--plane fields~\cite{El1}. Tight contact structures are much more
misterious, and difficult to classify. A contact structure on $Y$
is~\emph{virtually overtwisted} if its pull--back to some finite cover
of $Y$ becomes overtwisted, while it is called~\emph{universally
tight} if its pull--back to the universal cover of $Y$ is tight.

A contact three--manifold $(Y,\xi)$ is \emph{symplectically fillable},
or simply \emph{fillable}, if there exists a compact symplectic
four--manifold $(W, \omega)$ such that (i) $\partial W=Y$ as oriented
manifolds (here $W$ is oriented by $\om\wedge\om$) and (ii)
$\omega\vert_{\xi}\not=0$ at every point of $Y$. $(Y,\xi)$ is
symplectically \emph{semi}--fillable if there exists a fillable
contact manifold $(N, \eta )$ such that $Y\subset N$ and
$\eta|_Y=\xi$. Semi--fillable contact structures are tight~\cite{El3,
Gr}. The converse is known to be false by work of Etnyre and Honda,
who recently found two examples of tight but not semi--fillable
contact three--manifolds~\cite{EH}. Nevertheless, all such examples
known at present are virtually overtwisted, so it is natural to wonder
whether every universally tight contact structure is symplectically
semi--fillable.

In this paper we study certain virtually overtwisted tight contact
structures discovered by Ko Honda. Denote by $Y_{g,n}$ the total space
of an oriented $S^1$--bundle over $\Si_g$ with Euler number $n$.
Honda gave a complete classification of the tight contact structures
on $Y_{g,n}$~\cite{H2}. The three--manifolds $Y_{g,n}$ carry
infinitely many tight contact structures up to diffeomorphism. The
hardest part of the classification involves two virtually overtwisted
contact structures $\xi_0$ and $\xi_1$, which exist only when $n\geq
2g$. Honda conjectured that $\xi_0$ and $\xi_1$ are not symplectically
semi--fillable~\cite{H2}. The main theorem of the present paper
extends our earlier results regarding these structures~\cite{LS},
establishing Honda's conjecture:
\begin{thm}\label{t:main}
For $n\geq 2g>0$, the tight contact structures $\xi_0$ and $\xi_1$ on
$Y_{g,n}$ are not symplectically semi--fillable.
\end{thm}
The proof of Theorem~\ref{t:main} consists of two steps. In the first
step, we derive a contact surgery presentation for $\xi_0$ and $\xi_1$
in the sense of~\cite{DG2}, and we use it to determine the homotopy
type of $\xi_0$ and $\xi_1$ considered as oriented two--plane
fields. This is done in Sections~\ref{s:surgery} and~\ref{s:homotopy}.

In the second step, using specific properties of the $Spin ^c$
structures $\t _{\xi _i}$ on $Y_{g,n}$ induced by $\xi_i$ ($i=0,1$)
we generalize a result of the first author~\cite{L} so it applies to
the situation at hand. Using this generalization together with an
analytic computation of Nicolaescu's~\cite{N}, we are able to
determine the possible homotopy types of a semi--fillable contact
structure inducing either $\t_{\xi_0}$ or $\t_{\xi_1}$. This is done
in Section~\ref{s:proof}.

Theorem~\ref{t:main} follows immediately from the fact that the two
sets of homotopy classes determined in the two steps above have empty
intersection.

\section{Contact surgery presentations for $\xi _0$ and $\xi _1$}
\label{s:surgery}

A smooth knot $K$ in a contact three--manifold $(Y,\xi)$ which is
everywhere tangent to $\xi$ is called~\emph{Legendrian}. The contact
structure $\xi$ naturally induces a framing of $K$ called
the~\emph{contact framing}.

Let $\Si_g$ be a closed, oriented surface of genus $g\geq 1$, and let
$\pi\co Y_{g,n}\to\Si_g$ denote an oriented circle bundle over $\Si_g$
with Euler number $n$. Let $\xi$ be a contact structure on $Y_{g,n}$
such that a fiber $f=\pi^{-1}(s)\subset Y_{g,n}$ ($s\in\Si_g$) is
Legendrian. We say that $f$ has \emph{twisting number $k$} if the
contact framing of $f$ is $k$ with respect to the framing determined
by the fibration $\pi$.  A contact structure on $Y_{g,n}$ is called
\emph{horizontal} if it is isotopic to a contact structure transverse
to the fibers of $\pi$. 

Let $\ze$ be a horizontal contact structure on $Y_{g,2g-2}$ such that
a fiber $f$ of the projection $\pi$ is Legendrian with twisting number
$-1$ (the existence of such a contact structure is well--known,
cf.~\cite{Gi2},~\S 1.D). Let $n\geq 2g$, and view the bundle
$Y_{g,n}\to\Si_g$ as obtained by performing a $-\frac
1{p+1}$--surgery, where $p=n-2g+1$, along the fiber $f$ of $\pi\co
Y_{g,2g-2}\to\Si_g$ with respect to the trivialization induced by the
fibration $\pi$. It was observed by Honda~(\cite{H2},~\S 5) that there
are two possible ways of extending $\ze$ from the complement of a
standard neighborhood of $f$ to a tight contact structure on
$Y_{g,n}$. This determines the contact structures $\xi_0$ and $\xi_1$.

The construction of $\xi_0$ and $\xi_1$ can be viewed as a particular
case of a more general construction. In fact, given a Legendrian knot
$K$ in a contact three--manifold $(Y,\xi)$ and a rational number
$r\in\Q$, it is possible to perform a~\emph{contact $r$--surgery}
along $K$ to obtain a new contact three--manifold
$(Y',\xi')$~\cite{DG1, DG2}. Here $Y'$ is the three-manifold obtained
by a smooth $r$--surgery along $K$ with respect to the contact
framing, while $\xi'$ is constructed by extending $\xi$ from the
complement of a standard neighborhood of $K$ to a tight contact
structure on the glued--up solid torus. Such extension exists once
$r\neq 0$. In general there are several ways to extend $\xi$, but up
to isotopy there is only one if $r=\frac 1k$, $k\in\Z$, and two if
$r=\frac{p}{p+1}$ and $p>1$, as follows from~\cite{DG1},
Propositions~3, 4 and 7. When $r=-1$ the corresponding contact surgery
coincides with Legendrian surgery~\cite{El2,Go,W}. A simple
computation using the fact that the fiber $f$ of $Y_{g,2g-2}$ has
twisting number $-1$ with respect to the contact structure $\ze$ shows
that $\{\xi_0, \xi_1\}$ can be defined as the set of contact
structures obtainable by contact $\frac {p}{p+1}$--surgery along $f$.

\begin{figure}[ht]
\begin{center}
\epsfig{file=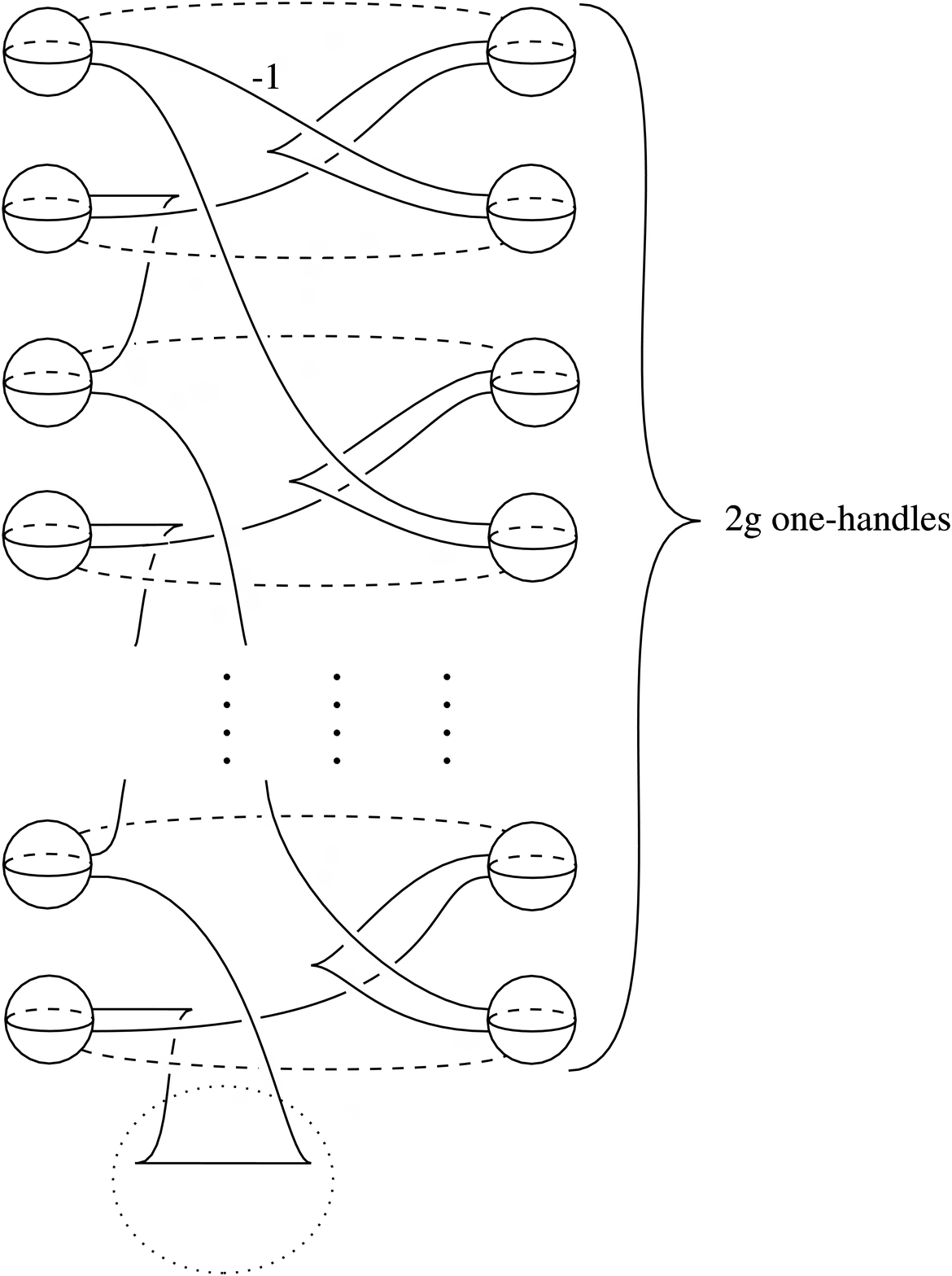, height=13cm}
\end{center}
\caption{The Stein four--manifold with boundary, 
diffeomorphic to $D_{g,2g-2}$}
\label{f:figure1}
\end{figure}

{From} now on, we shall indicate a contact $r$--surgery along a
Legendrian knot $K$ by writing the coefficient $r$ next to it.
Consider the result of performing contact $(-1)$--surgery on the
Legendrian knot in standard form in Figure~\ref{f:figure1} (here we
are using the notation of~\cite{Go}, see especially Definition~2.1).
Since contact $(-1)$--surgery is equivalent to Legendrian surgery,
Figure~\ref{f:figure1} also represents a Stein four--manifold $W$ with
boundary~\cite{Go}. As a smooth four--manifold, $W$ is diffeomorphic
to the two--disc bundle $D_{g,2g-2}$ with Euler number $2g-2$ over a
surface of genus $g$. This can be checked by converting the contact
surgery coefficient into the corresponding smooth surgery coefficient
e.g.~via the formulas found in~\cite{Go} or~\cite{GS}. Since by
construction the boundaries of Stein four--manifolds come equipped
with Stein fillable contact structures, we have a Stein fillable
contact structure $\ze(g)$ on $Y_{g,2g-2}$, which is tight
by~\cite{El3, Gr}.

\begin{lem}\label{l:horizontal}
The contact structure $\ze(g)$ is horizontal. Moreover, after an
isotopy the map $\pi\co Y_{g,2g-2}\to\Si_g$ has a fiber with twisting
number equal to $-1$.
\end{lem}

\begin{proof}
The existence of a Legendrian knot isotopic to a fiber with twisting
number $-1$ is apparent from Figure~\ref{f:figure1}. On the other
hand, contact $(-1)$--surgery on a Legendrian knot isotopic to a fiber
and having twisting number $\geq 0$ would result in a Stein manifold
containing a sphere with self-intersection $\geq -1$, contradicting
the adjunction inequality for Stein manifolds~\cite{LM}. By the
classification of tight contact structures on $Y_{g,2g-2}$
with~\emph{negative twisting number} i.e.~such that the twisting number of
any closed Legendrian curve isotopic to a fiber is $<0$ (\cite{H2},
Theorem~2.11), we conclude that the diagram of Figure~\ref{f:figure1}
represents a horizontal contact structure.
\end{proof}

\begin{figure}[ht]
\begin{center}
\epsfig{file=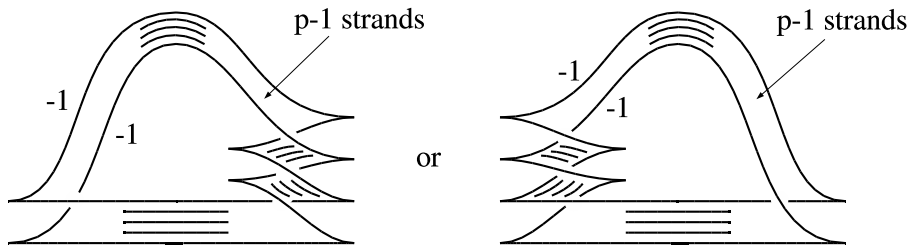, height=3.6cm}
\end{center}
\caption{The Legendrian surgeries equivalent to contact
$-\frac{p}{p-1}$--surgery}
\label{f:figure2}
\end{figure}

By~\cite{DG2}, Proposition~3, any contact $r$--surgery with $r<0$ is
equivalent to a Legendrian surgery along a Legendrian link. Moreover,
the set of Legendrian links which correspond to some contact
$r$--surgery is determined via a simple algorithm by the Legendrian
knot and the continued fraction expansion of $r$. For example, let $K$
be a Legendrian unknot in the standard contact three--sphere with
Thurston--Bennequin invariant equal to $-1$. Then, a contact
$-\frac{p}{p-1}$--surgery ($p>1$) along $K$ is equivalent to
Legendrian surgery along one of the Legendrian links in 
Figure~\ref{f:figure2}.
According to~\cite{DG2}, Proposition~7, a contact
$\frac{p}{p+1}$--surgery on a Legendrian knot $K$ is equivalent to a
contact $\frac{1}{2}$--surgery on $K$ followed by a contact
$-\frac{p}{p-1}$--surgery on a Legendrian push--off of $K$.
By~\cite{DG1}, Proposition~9, a contact $\frac{1}{2}$--surgery on a
Legendrian knot $K$ can be replaced by two contact $(+1)$--surgeries,
one on $K$ and the other on a Legendrian push--off of $K$.

\begin{figure}[ht]
\begin{center}
\epsfig{file=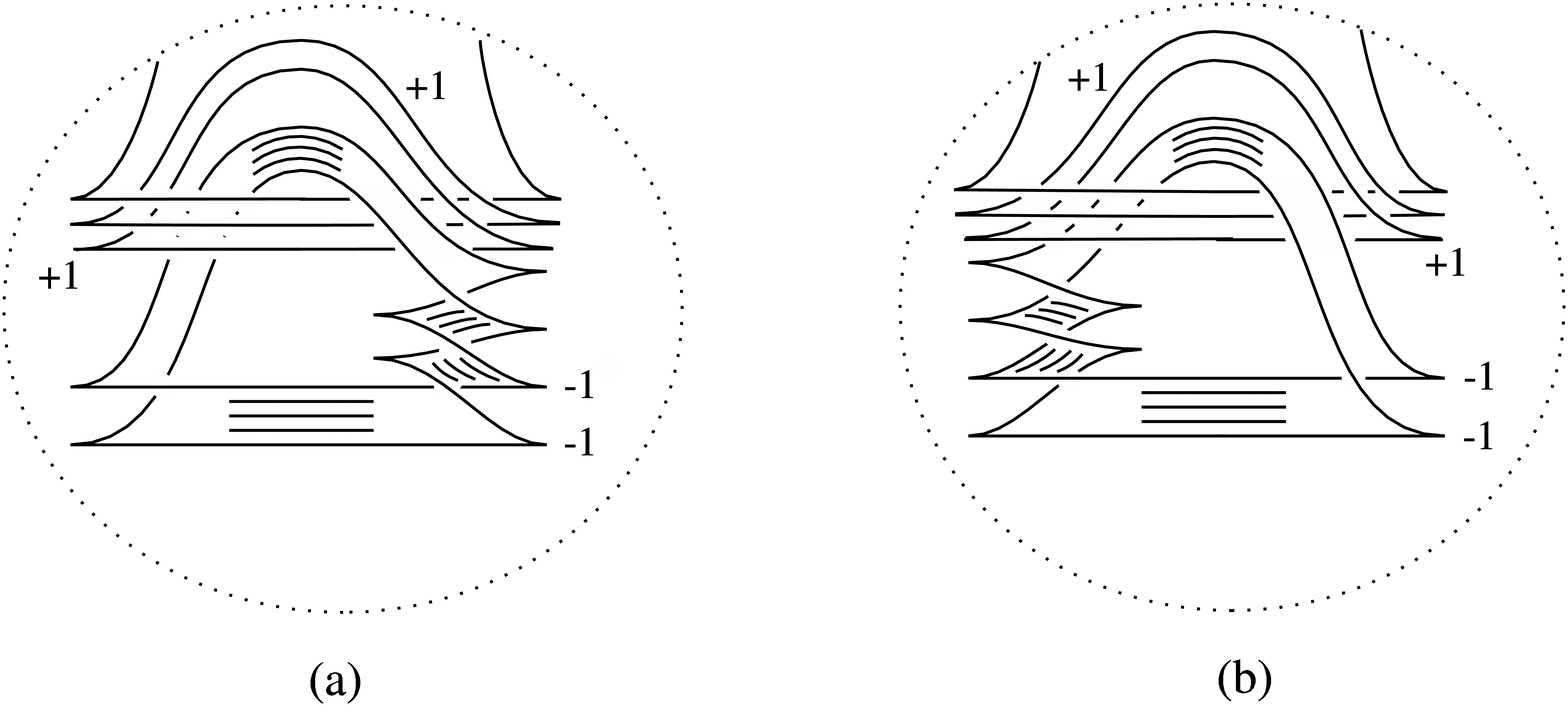, height=6cm}
\end{center}
\caption{Pictures to be pasted in Figure~\ref{f:figure1} to obtain
$\xi_0$ and $\xi_1$}
\label{f:figure3}
\end{figure}

This implies that if we perform a Legendrian $\frac{p}{p+1}$--surgery
on a Legendrian fiber of $(Y_{g,2g-2},\ze(g))$ with twisting number
$-1$, the resulting contact structures will have contact surgery
presentations obtained by replacing the ``dotted ellipse'' in
Figure~\ref{f:figure1} with either Figure~\ref{f:figure3}(a)
or~\ref{f:figure3}(b). More precisely, we can define $\xi_0$, respectively
$\xi_1$, as the contact structure obtained by using
Figure~\ref{f:figure3}(a), respectively Figure~\ref{f:figure3}(b).

\section{Homotopy classes of $\xi_0$ and $\xi_1$}
\label{s:homotopy}

\subsection*{Homotopy theory of oriented two--plane fields 
on three--manifolds} Let $\Xi_Y$ denote the space of oriented
two--plane fields on the closed, oriented three--manifold $Y$. Since a
$Spin^c$ structure on a three--manifold can be interpreted as an
equivalence class of nowhere vanishing vector fields~\cite{Tu}, by
taking the oriented normal, a two--plane field $\xi\in\Xi_Y$ naturally
induces a $Spin^c$ structure $\t_\xi$, which depends only on the homotopy
class $[\xi]$.  Therefore there is a map $p\co\pi_0(\Xi_Y)\to
Spin^c(Y)$ defined as $p([\xi])=\t_\xi$. It is not difficult to show
that, if $Y$ is connected, there is a non--canonical identification of
each fiber $p^{-1}(\t_\xi)$ with $\Z/d(\t_\xi)\Z$, where
$d(\t_\xi)\in\Z$ is the divisibility of $c_1(\xi)\in H^2(Y; \Z )$, and
is zero if $c_1(\xi)$ is a torsion element (see, e.g.~\cite{Go},
Proposition~4.1).

When $c_1(\xi)$ is torsion the two--plane fields inducing the same
$Spin^c$ structure $\t_{\xi}$ can be distinguished by a numerical
invariant. Suppose that $X$ is a compact 4-manifold with $\partial X =
Y$, with $X$ carrying an almost--complex structure $J$ whose complex
tangents at the boundary form an oriented two--plane field
homotopic to $\xi$ on $Y$. Observe that the fact that $c_1(\xi)$ is
torsion implies that $c^2_1(X,J)\in\Q$ makes sense. 
\begin{thm}[\cite{Go}]
The rational number
\[
d_3(\xi) = \frac{1}{4}(c^2_1(X, J) - 3\sigma
(X)-2\chi (X))\in \Q
\] 
depends only on $[\xi]$, not on the almost--complex four--manifold
$(X, J)$. Moreover, two two--plane fields $\xi_1$ and $\xi_2$ inducing
the same $Spin^c$ structure with torsion first Chern class are
homotopic if and only if $d_3(\xi _1)=d_3(\xi _2)$. \qed
\end{thm}

In the following we shall refer to the invariant $d_3$ as
the~\emph{three--dimensional invariant}.

\subsection*{Attaching two--handles and homotopy invariants}
Recall that contact $(-1)$--surgery, i.e. Legendrian surgery, can be
viewed as the result of attaching a symplectic two--handle~\cite{W}.
In fact, attaching the two--handle to a contact three--manifold
$(Y_1,\xi_1)$ gives rise to a cobordism $W$ between $Y_1$ and the
three--manifold underlying the contact three--manifold $(Y_2,\xi_2)$
resulting from the three--dimensional contact surgery. Furthermore,
$W$ carries an almost--complex structure whose complex tangent lines at
the boundary coincide with $\xi_1$ and $\xi_2$ (see e.g.~\cite{EH2}). 

In the case of contact $(+1)$-surgery, there is still a smooth
cobordism $W$ between $Y_1$ and $Y_2$. One can easily check the
existence of an almost--complex structure $J$ on the complement of a
ball $B$ in the interior of $W$, with $J$ inducing $\xi_1$ and $\xi_2$
as tangent complex lines. We define $q$ to be the three--dimensional
invariant of the two--plane field induced by $J$ on $\partial
B$. Observe that, although $J$ may not extend to the whole cobordism,
$J$ induces a $Spin^c$ structure $\s_J$ which does extend -- uniquely
-- to $W$.
\begin{lem}\label{l:value}
The value of $q$ is $\frac 12$.
\end{lem}

\begin{proof}
Consider an oriented Legendrian unknot $K$ in the standard contact
three--sphere with Thurston--Bennequin invariant equal to $-1$ and
vanishing rotation number. We view the standard contact three--sphere
as the contact boundary of the unit ball $B_1(0)\subset\C^2$. Attach a
smooth two--handle $H_1$ to $B_1(0)$ with framing $+1$ with respect to
the contact framing. The result is a smooth four--manifold $X$
diffeomorphic to $S^2\x D^2$. The unique $Spin^c$ structure on
$B_1(0)$ extends to a $Spin^c$ structure $\s$ on $X$, restricting to
$H_1$ as the $Spin^c$ structure defined above. Denote by $k$ the
value of $c_1(\s)$ on a generator of the second homology
group of $X$.

Let $K'$ be a Legendrian push--off of $K$, which we may assume
disjoint from $H_1$, and attach a symplectic two--handle $H_2$ to $K'$
realizing Legendrian surgery on $K'$. The $Spin^c$ structure $\s$
extends over $H_2$, and the value of its first Chern class on the
homology generator corresponding to $K'$ is $0$, because $K'$ has
vanishing rotation number (see~\cite{Go}, especially the proof of
Proposition~2.3). By~\cite{DG1}, Proposition~8, the resulting contact
three--manifold is just the standard contact three--sphere. Its
three--dimensional invariant $d_3$ is $-\frac 12$, but when viewed as the
result of the above construction, $d_3$ can also be expressed as
$\frac{1}{4}(2k^2-4)+q$.

We can generalize this argument using Legendrian push--offs
$K_1$, $K_1', \ldots$, $K_n$, $K_n'$ of $K$ by performing contact
$(+1)$--surgeries on $K_1,\ldots , K_n$ and contact $(-1)$--surgeries
on $K_1', \ldots , K_n'$. The resulting contact three--manifold is the
standard contact three--sphere again. A homological computation as 
before gives the identity
\[
\frac{1}{4}(n+1)(k^2n-2)+nq=-\frac 12,
\]
which must hold for all $n\in\N$. This implies that $k=0$ and
$q=\frac 12$.
\end{proof}

\subsection*{$Spin^c$ structures on disc and circle bundles}
Let $D_{g,n}$ be the oriented disc bundle with Euler number $n$ over a
closed oriented surface of genus $g$. By e.g. fixing a metric on
$D_{g,n}$ one sees that the tangent bundle of $D_{g,n}$ is isomorphic
to the direct sum of the pull--back of $T\Si_g$ and the vertical
tangent bundle, which is isomorphic to the pull--back of the real
oriented two--plane bundle $E_{g,n}\to\Si_g$ with Euler number $n$. In
short, we have
\begin{equation}\label{e:split}
TD_{g,n}\cong \pi^*(T\Si_g\oplus E_{g,n}).
\end{equation}
This splitting of $TD_{g,n}$ naturally endows $D_{g,n}$ with and
almost--complex structure which induces a $Spin^c$ structure $\s_0$
on $D_{g,n}$.  The orientation on $D_{g,n}$ determines an isomorphism
$H^2(D_{g,n};\Z)\cong\Z$, so the set $Spin^c(D_{g,n})=\s_0 +
H^2(D_{g,n};\Z)$ can be canonically identified with the integers. We
denote by $\s_e=\s_0+e\in Spin^c(D_{g,n})$ the element corresponding
to the integer $e\in\Z\cong H^2(D_{g,n};\Z)$.

Consider $Y_{g,n}=\del D_{g,n}$. We have $H_1(Y_{g,n};\Z)\cong
H^2(Y_{g,n}; \Z) \cong \Z ^{2g}\oplus \Z/n\Z$, where the summand
$\Z/n\Z$ is generated by the Poincar\'e dual $F$ of the class of a
fiber of the projection $\pi\co Y_{g,n}\to\Si_g$. Each $Spin^c$
structure $\s_e\in Spin^c(D_{g,n})$ determines by restriction a
$Spin^c$ structure $\t_e\in Spin^c(Y_{g,n})$ with $\t_e=\t_0+eF$,
$e\in\Z$. Since $nF=0$, we see that $\t_{e+n}=\t_e$ for every
$e$. Therefore, $\t_0,\ldots,\t_{n-1}$ is a complete list of
\emph{torsion $Spin^c$ structures} on $Y_{g,n}$, i.e. $Spin^c$
structures on $Y_{g,n}$ with torsion first Chern class. In short, the
$Spin^c$ structures on $Y_{g,n}$ which extend to the disc bundle are
precisely the torsion ones.

\subsection*{Homotopy invariants of the contact structures $\xi _i$}
\begin{figure}[t]
\begin{center}
\epsfig{file=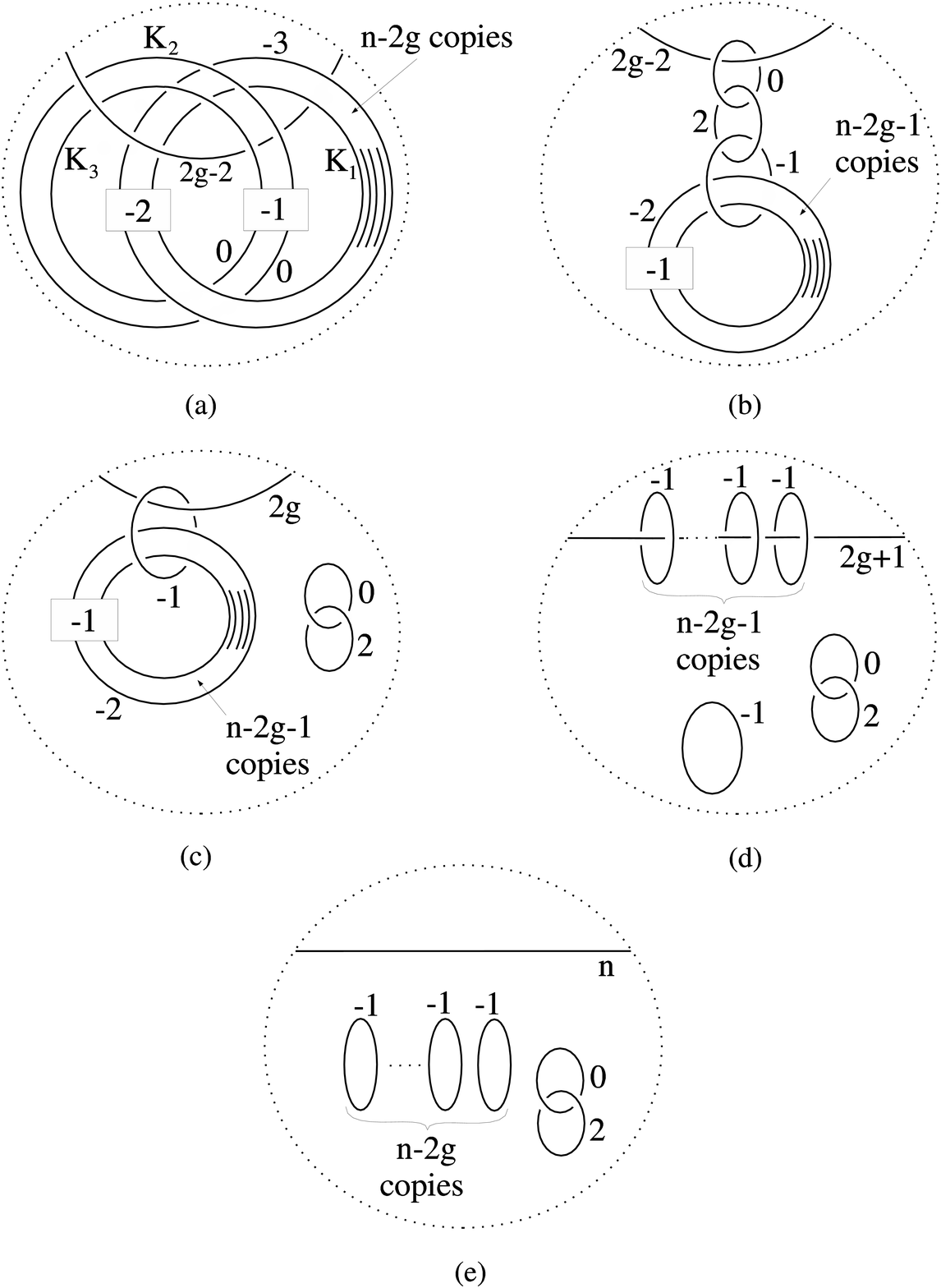, height=16cm}
\end{center}
\caption{The diffeomorphism between $X$ and $D_{g,n}\# S^2\x S^2
\# (n-2g){\overline\CP}^2$}
\label{f:figure4}
\end{figure}

Let $W$ be the Stein four--manifold with boundary diffeomorphic to
$D_{g, 2g-2}$ as given by Figure~\ref{f:figure1}. Consider the smooth
four--dimensional handlebody $X$ obtained by attaching to $W$ the
two--handles realizing the contact surgeries described in
Figure~\ref{f:figure3}(a) or ~\ref{f:figure3}(b). Converting the
contact framing coefficients into the usual ones, we see that a framed
link presentation of $X$ is obtained by pasting
Figure~\ref{f:figure4}(a) in place of the `dotted ellipse' in
Figure~\ref{f:figure1}.

By the discussion above on attaching two--handles we know that,
corresponding to each of Figure~\ref{f:figure3}(a) and
~\ref{f:figure3}(b), there is an almost--complex structure on $X$ minus
two balls lying in the interior of the two--handles realizing the
$(+1)$--surgeries. Moreover, the two almost--complex structures
determine the two--plane fields $\xi_0$ and $\xi_1$ on $\del X$ and
two $Spin^c$ structures $\s_0$ and $\s_1$ on $X$. Observe that, since
the rotation number of the Legendrian knot in Figure~\ref{f:figure1}
vanishes, it follows from~\cite{Go}, Theorem~4.12, that $c_1(W)=0$. In
the same way, it follows that we can choose an orientation of the
$n-2g$ linking knots with framing $-3$ in Figure~\ref{f:figure4}(a) so
that $c_1(\s_i)$ evaluates as $(-1)^i$ on all the corresponding
homology classes. Finally, by the argument given in the proof of
Lemma~\ref{l:value}, $c_1(\s_i)$ evaluates trivially on the generators
of $H_2(X;\Z)$ determined by the two--handles corresponding to the
$(+1)$--surgeries.

The four--manifold $X$ is diffeomorphic to $D_{g,n}\# S^2\x S^2\#
(n-2g){\overline\CP}^2$. One can see this by performing a sequence of
handleslides on the Kirby diagram as shown in Figure~\ref{f:figure4}.
In fact, start by sliding over the knot $K_1$ in
Figure~\ref{f:figure4}(a) the remaining $(n-2g-1)$ $(-3)$--framed
circles. Then, slide $K_1$ over $K_2$ and finally $K_2$ over $K_3$,
obtaining~\ref{f:figure4}(b). Sliding the long $(2g-2)$--framed arc
over the $2$--framed knot and using the $0$--framed normal circle to
separate the $2$--framed circle from the rest of the diagram, we
get~\ref{f:figure4}(c). Blowing down the $(-1)$--circle results
in~\ref{f:figure4}(d), and $(n-2g-1)$ further blow downs
give~\ref{f:figure4}(e).
Following the handle slides of Figure~\ref{f:figure4} on the
homological level we see that $c_1(\s_i)$ evaluates on the generator of
the second homology of $D_{g,n}$ as $(-1)^i(n-2g)$. Moreover, it
evaluates as $(-1)^i$ on generators of the ${\overline\CP}^2$
summands, and vanishes when restricted to the $S^2\x S^2$
summand. This immediately implies that the $Spin^c$ structure
$\t_{\xi_i}$ is equal to the restriction of the unique $Spin^c$
structure $\s_e\in Spin^c(D_{g,n})$ such that $c_1(\s_e)$ evaluates on
the generator of $H_2(D_{g,n};\Z)$ as $(-1)^i(n-2g)$.  Since
the value of $c_1(\s_0)$ on the generator is $2-2g+n$, $e$ satisfies
the equation:
\[
2-2g+n+2e=(-1)^i (n-2g).
\]
Therefore we get $e=-1$ or $e=2g-1+n$ respectively for $i=0$ or
$i=1$. Since $\s_e|_{Y_{g,n}}=\t_e$, we conclude that
\[
\t_{\xi _i }=\t_{2ig-1}
\]
for $i=0,1$. Observe that this result is consistent with the
independent calculation made in~\cite{LS}.

\begin{lem}\label{l:d3ofxi}
The value of the three--dimensional invariant of $\xi_i$ is
\[
d_3(\xi_i)=\frac{n^2-3n+4g^2}{4n}.
\]
\end{lem}
\begin{proof}
We have $\chi(X)=n-4g+4$ and $\sigma (X)=1-n+2g$. {From} what we know about
$c_1(\s_i)$ is easy to deduce that
\[
c_1^2(\s_i)=-2\frac{g(n-2g)}{n}.
\]
In order to compute the three--dimensional invariant we need to take
into account the correction term $q$ for each of the two contact
$(+1)$--surgeries. Using Lemma~\ref{l:value} we conclude
\[
d_3(\xi _i)
=\frac{1}{4}(c_1^2(\s_i)-2\chi(X)-3\si(X))+2=
\frac{n^2-3n+4g^2}{4n}. 
\]
\end{proof} 

\section{The proof of Theorem~\ref{t:main}}
\label{s:proof}

\begin{thm}\label{t:pani}
Let $n\geq 2g>0$, and let $\xi$ be a two--plane field on $Y_{g,n}$
such that $\t_{\xi}\in\{\t_{\xi_0},\t_{\xi_1}\}$. If $\xi$ is homotopic
to a semi--fillable contact structure, then 
\[
d_3(\xi )=\frac{n^2+n+4g^2}{4n}-2g-2.
\]
\end{thm}

\begin{proof}
In the proof of Theorem~2.1 of~\cite{L} it is shown that if $Y$ is 
a closed three--manifold and $\t\in Spin^c(Y)$ is torsion and satisfies:
\begin{itemize}
\item all Seiberg-Witten solutions in $\t$ are reducible and
\item the moduli space of the Seiberg-Witten solutions in $\t$ is a smooth
manifold and the corresponding Dirac operators have trivial kernels,
\end{itemize}
then the expected dimension $d_1$ of the Seiberg-Witten moduli
space of solutions over a potential symplectic semi--filling of
$(Y_{g,n}, \xi_i)$ equipped with a cylindrical end metric and fixed
asymptotic limit is equal to $-1-b_1(Y_{g,n})$.

The moduli space of Seiberg-Witten solutions on $Y_{g,n}$ has been
determined in \cite{MOY} (see also \cite{OSz}). These results show
that the assumptions listed above hold for the moduli spaces
associated to the $Spin ^c$ structures $\t_{\xi_i}$. Therefore, the
conclusion $d_1=-1-b_1(Y_{g,n})$ holds. This implies that for each
$i=0,1$, $Y_{g,n}$ carries only one homotopy type of two--plane field
which contains potentially semi--fillable contact structures inducing
$\t_{\xi_i}$, because $d_1$ is equal to the three--dimensional
invariant plus an expression involving some topological terms and an
$\eta$--invariant~(\cite{L0}, Formula~3.1). In fact, such an
expression has been explicitely calculated in~\cite{N}, in the formula
preceding~(3.29), so our proof reduces to translating that formula
into our notations.

In Nicolaescu's notations the integer $\kappa $ corresponds to
$\t_{g-1+\kappa}$. This is because his ``base'' $Spin^c$ structure is
induced by a $Spin$ structure on $Y_{g,n}$ with associated bundle of
spinors ${\mathbb S}=\pi^* K_{\Si_g}^{-\frac 12}\oplus\pi^*
K_{\Si_g}^{\frac 12}\to Y_{g,n}$ (see text following Formula (2.6)
in~\cite{N}), and $\mathbb S$ is the restriction of $TD_{g,n}\otimes
\pi^* K_{\Si_g}^{\frac 12}\to D_{g,n}$ to the boundary.

The result we need is obtained by substituting $n$ for $\ell$ and $g$
or $n-g$ in place of $\kappa$ into the formula preceding~(3.29)
of~\cite{N}. (The formula we are using here differs from
Formula~(3.29) by the additive term $2g-1$ because (3.29) computes the
dimension of the whole moduli space rather than the dimension of the
moduli space of solutions with a fixed asymptotic limit,
i.e.~$d_1$). Explicitely, in our notation we have:
\[
-1-b_1(Y_{g,n})=d_1=d_3(\xi ) -\frac{1}{2}(2g-1)-\frac{1}{4}(n-1)-
\frac{\kappa ^2}{n} +\kappa
\]
where $b_1(Y_{g,n})=2g$ and the value of $\kappa $ to be substituted
is either $g$ or $n-g$ according to whether $\t_\xi=\t_{\xi_1}$ or
$\t_\xi=\t_{\xi_0}$, respectively. In both cases we obtain for
$d_3(\xi)$ the value given in the statement.
\end{proof}

\begin{proof}[Proof of Theorem~\ref{t:main}]
Let $\xi$ be a two--plane field representing a homotopy class inducing
$\t_{\xi _i}$ which might be represented by a semi--fillable contact
structure. Then, by Theorems~\ref{l:d3ofxi} and~\ref{t:pani} we have
$d_3(\xi_i) - d_3 (\xi)=2g+1>0$. Therefore, the homotopy classes
$[\xi_i]$ cannot be represented by semi--fillable contact structures.
\end{proof}

\noindent
{\bf Remarks.}
\noindent
(1) For $n<2g$ the circle bundle $Y_{g,n}$ admits no $Spin ^c$
structure for which the Seiberg-Witten moduli space has the properties
required by the proof of Theorem~\ref{t:pani}.

(2) The assumption $g>0$ in Theorem~\ref{t:pani} is necessary ---
$Y_{0,n}$ is a lens space on which all tight contact structures are
Stein fillable. The proof of Theorem~\ref{t:pani} breaks down since
the formula from~\cite{N} used in the proof holds only for $g\geq 1$.

(3) Notice that for $n=2g$ the two contact structures $\xi_0$ and $\xi _1$
coincide.

\end{document}